\documentclass{amsart}

\usepackage{xcolor}
\usepackage{mathtools}
\usepackage{amsmath}
\usepackage{amsthm}
\usepackage{graphicx}
\usepackage{amssymb}
\usepackage{epstopdf}
\usepackage{nicefrac}
\usepackage{mathrsfs}
\usepackage{hyperref}
\usepackage{enumitem}
\usepackage{tikz}
\usepackage{caption}
\usepackage{subcaption}
\usepackage{cite}
\usetikzlibrary{calc,arrows.meta,decorations.pathreplacing}

\theoremstyle{plain}
\newtheorem{theo}{Theorem}

\newtheorem*{q}{Question}

\newtheorem*{ack}{Acknowledgements}
\newtheorem*{ai}{Use of Artificial Intelligence}

\newtheorem*{example*}{Example}

\DeclareMathOperator{\Ha}{H}

\DeclareMathOperator{\spans}{span}

\title[Failure of Nehari for Multiplicative Hankel]{Failure of Nehari's Theorem for Multiplicative Hankel Forms in Schatten Classes II}
\author[K. Bampouras]{Konstantinos Bampouras}
\address{Department of Mathematics and Systems Analysis, School of Science, Aalto University, Espoo, Finland}
\email{konstantinos.bampouras@aalto.fi}

\date{\today}

\begin{document}
	
	\begin{abstract}
		We prove that for any $p>4$ we can find a Schatten class Hankel operator $\Ha_\phi$ defined on the Hardy space of the infinite torus $H^2(\mathbb{T}^\infty)$ that belongs to the Schatten class $S^p$ but does not admit a bounded symbol, extending the previously known range $p>(1-\log\pi /\log 4)^{-1}=5.7388...$.
	\end{abstract}
	
	\maketitle

	\section{Introduction}
For a sequence $\varrho=(\varrho_n)_{n\geq 1}\in \ell^2$, the associated multiplicative Hankel form is
$$\varrho(a,b)=\sum_{n,m=1}^\infty \varrho_{nm}a_n\overline{b_m},$$
initially defined for finitely supported $a,b\in\ell^2$. One may ask which properties of the sequence $(\varrho_n)_{n\geq1}$ determine boundedness of this form or membership of the associated operator in a Schatten class. 

For $k\in\mathbb N$, write its prime factorization as
$$k=p_1^{\alpha_1}p_2^{\alpha_2}...,$$ where $(p_n)_{n\geq 1}$ is the sequence of primes, $(\alpha_n)_{n\geq 1}$ is a finite sequence of nonnegative integers,
and set
$$\alpha(k)=(\alpha_1,\alpha_2,...).$$
For $x=(x_n)_{n\geq1}\in\ell^2$, define $\phi_x\in H^2(\mathbb T^\infty)$ (see the next paragraph for the precise definition of this Hardy space) by
$$\phi_x(z)=\sum_{n=1}^\infty x_n z^{\alpha(n)}.$$
Since the monomials $z^{\alpha(n)}$ form an orthonormal basis, we have
$$\|\phi_x\|_{H^2(\mathbb T^\infty)}=\|x\|_{\ell^2},\qquad\widehat{\phi}_x(\alpha(n))=x_n.$$
For a function $\phi\in H^2(\mathbb{T}^\infty)$, its corresponding Hankel operator $\Ha_\phi:H^2(\mathbb{T}^\infty)\to H^2(\mathbb{T}^\infty)$ is given by the matrix $(\widehat{\phi}(\alpha(n)+\alpha(m)))_{n,m\geq 1}$ acting on the base $(z^{\alpha(n)})_{n\geq 1}$, initially defined for functions with finitely supported Fourier transform.
Therefore, for $\varrho\in\ell^2$ and finitely supported $a,b\in\ell^2$,
$$
	\left\langle \Ha_{\phi_\varrho}\phi_a,\phi_b\right\rangle
	=
	\sum_{n,m=1}^\infty
	\widehat{\phi}_\varrho\bigl(\alpha(n)+\alpha(m)\bigr)
	\widehat{\phi}_a\bigl(\alpha(n)\bigr)
	\overline{\widehat{\phi}_b\bigl(\alpha(m)\bigr)}\\
	=
	\sum_{n,m=1}^\infty
	\varrho_{nm}a_n\overline{b_m},
$$
where in the last equality we used $\alpha(n)+\alpha(m)=\alpha(nm)$. Thus, the multiplicative Hankel form associated with $\varrho$ is identified with the Hankel form associated with $\phi_\varrho$ on $H^2(\mathbb T^\infty)$. We will therefore work in this setting. For background on this correspondence, we refer the reader to \cite{helson2005dirichlet,queffelec2013diophantine}.

	On the infinite torus $\mathbb{T}^\infty=\mathbb{T}\times \mathbb{T}\times...$ we consider the Haar measure $m_\infty$, which is the countably infinite product of the normalised Lebesgue measure $m_1$ on $\mathbb{T}$. The Fourier coefficients of an $L^1(\mathbb{T}^\infty)$ function are defined on the set $\mathbb{Z}^{(\infty)}$ of all finitely supported sequences of integers, by
	$$\hat{f}(\alpha)=\int_{\mathbb{T}^\infty} f(z)z^{-\alpha} dm_\infty(z),$$ where $\alpha=(\alpha_1,...,\alpha_k,0,...)\in \mathbb{Z}^{(\infty)}$ and $z^{-\alpha}=z_1^{-\alpha_1}...z_k^{-\alpha_k}.$
	The Hardy spaces of the infinite torus, for $p\geq 1$, are defined by
	$$H^p(\mathbb{T}^\infty)=\{f\in L^p(\mathbb{T}^\infty):\hat{f}(\alpha)=0 \text{ for all }\alpha\notin \mathbb{N}_0^{(\infty)}\},$$
	endowed with the $L^p(\mathbb{T}^\infty)$ norm.
	For a symbol $\phi\in L^2(\mathbb{T}^\infty)$, the associated Hankel operator $\Ha_\phi:H^2(\mathbb{T}^\infty)\to H^2(\mathbb{T}^\infty)$ is initially defined on analytic polynomials (which are dense) by
	$$\widehat{\Ha_\phi f}(\alpha) = \sum_{\beta\in \mathbb{N}_0^{(\infty)}}\widehat{\phi}(\alpha+\beta)\hat{f}(\beta), \quad\alpha\in \mathbb{N}_0^{(\infty)}.$$ Analogously, we can define Hardy spaces and Hankel operators in finite dimensions and regard $H^p(\mathbb{T}^n)$ as a subspace of $H^p(\mathbb{T}^\infty)$. For a function $\phi\in L^\infty(\mathbb{T}^\infty)$, it is evident that
	$$\|\Ha_\phi\|\leq \|\phi\|_{L^\infty}.$$ Therefore, the very natural question arises: is it true that for every bounded Hankel operator $\Ha_\phi$, we can find $\psi\in L^\infty(\mathbb{T}^\infty)$ such that $\Ha_\phi=\Ha_\psi$? 
	This question is known as Nehari's problem, since Nehari gave a positive answer for the case of dimension $1$ \cite{MR82945}. The higher dimensional cases are open. If we assume that Nehari's problem holds, then by the open mapping theorem, we can find $C>0$ such that for every $f\in H^1(\mathbb{T}^\infty)$ and $\phi\in L^2(\mathbb{T}^\infty)$ with $\Ha_\phi$ bounded, 
	$$|\langle f,\phi\rangle|\leq C\|\Ha_\phi\|\|f\|_{L^1}.$$
	Since $H^p(\mathbb{T}^n)$ is a subspace of $H^p(\mathbb{T}^\infty)$, the Nehari problem would also have a positive answer for every finite dimension, and we could find constants $C_n\leq C$ with 
	$$|\langle f,\phi\rangle|\leq C_n\|\Ha_\phi\|\|f\|_{L^1}.$$
	Ortega-Cerd\`{a} and Seip \cite{MR2993031} showed that $$C_{2n}\geq \left(\frac{\pi^2}{8}\right)^{n/2},$$ disproving Nehari's theorem on the infinite torus. Brevig \cite{Brevig2022MinimalNorm} also strengthened the bound to 
	$$C_{2n}\geq \left(\frac{5\pi}{\pi+6\sqrt 3}\right)^n.$$
	
	Instead of asking whether any bounded Hankel operator admits a bounded symbol, Helson \cite{MR2263964} asked whether this consequence holds for every Hilbert-Schmidt Hankel operator and gave a positive answer. His proof heavily relied on the fact that we can compute the Hilbert-Schmidt norm to be
	$$\|\Ha_\phi\|_{S^2}^2=\sum_{\alpha\in\mathbb{N}_0^{(\infty)}}|\widehat{\phi}(\alpha)|^2d(\alpha),$$ where $$d(\alpha_1,...,\alpha_k)=(\alpha_1+1)...(\alpha_k+1),$$ and thus a duality argument gives that Nehari's consequence for Hilbert-Schmidt operators is equivalent to the inequality
	\begin{equation}\label{eq1} \sum_{\alpha\in\mathbb{N}_0^{(\infty)}}\frac{|\hat{f}(\alpha)|^2}{d(\alpha)}\leq C\|f\|_{L^1}^2,\quad f\in H^1(\mathbb{T}^\infty).
	\end{equation}
	Helson used the one-dimensional Carleman's inequality
	$$\sum_{k\geq 0}\frac{|\hat{f}(k)|^2}{k+1}\leq \|f\|_{L^1}^2,\quad f\in H^1(\mathbb{T}),$$ and induction to prove inequality \eqref{eq1} with $C=1$.
	
	In 2015, Brevig and Perfekt \cite{BrevigPerfekt2015} asked the following more general question
	\begin{q}
		For which $p>2$ is it true that every Schatten class Hankel operator $\Ha_\phi$ in $S^p(H^2(\mathbb{T}^\infty))$ admits a bounded symbol?
	\end{q}
	Other than Helson's inequality, no positive result is known in dimensions greater than $2$. In dimension $2$, however, the recent result \cite{bampouras2026recoveringproductbmoschatten} gives a positive answer for every $p>2$.

	In the same paper, Brevig and Perfekt \cite{BrevigPerfekt2015} gave a negative answer for every $p>(1-\log \pi/\log 4)^{-1}=5.7388...$. In this paper, we are going to extend their result by proving the following theorem.
	\begin{theo}\label{thm:main}
		For every $p>4$ there exists a Hankel operator $\Ha_\phi\in S^p(H^2(\mathbb{T}^\infty))$ that does not admit a bounded symbol.
	\end{theo}
	Before concluding the introduction, we mention a related result. It was recently shown that the Nehari consequence also fails for all $p>4$ for Hankel operators on Paley-Wiener spaces of sufficiently smooth sets \cite{bampouras2026nehari}. It would be interesting to understand whether the appearance of the same threshold in these two settings reflects some connection between them. In particular, the endpoint case $p=4$ seems especially interesting, since its Schatten norm can be explicitly computed.
	\begin{ack}
		The author was supported by the Research Council of Finland, project no. 371637. They are a member of the Finnish Centre of Excellence in Randomness and Structures. The author would also like to thank Viktor Andersson, Ole Fredrik Brevig and Karl-Mikael Perfekt for their helpful comments regarding technique and presentation.
	\end{ack}
	\begin{ai}
		ChatGPT (OpenAI) suggested introducing the auxiliary function \(D_\eta\) and using its weighted integral \(I_\kappa(\eta)\) in the proof. It also assisted with intermediate estimates and preliminary proof drafts. The author checked and revised the arguments, simplified the exposition, and wrote the final manuscript. The author takes full responsibility for the mathematical content and any remaining errors.
	\end{ai}
	\section{Proof of Theorem \ref{thm:main}}
	For contradiction, let us assume that for some $p>4$ every Hankel operator $\Ha_\phi\in S^p(H^2(\mathbb{T}^\infty))$ admits a bounded symbol. Then, as proved in \cite{BrevigPerfekt2015} by Brevig and Perfekt, for every $f\in H^1(\mathbb{T}^n)$ and every symbol $\phi\in H^2(\mathbb{T}^n)$ with $\Ha_\phi\in S^p(H^2(\mathbb{T}^n))$ it must hold that
	\begin{equation}\label{eq2} |\langle f,\phi\rangle|\leq \|\Ha_\phi\|_{S^p}\|f\|_{L^1(\mathbb{T}^n)}.
	\end{equation}
	We will denote by $m_n$ the product measure on $\mathbb{T}^n$, that is, the product of the normalised Lebesgue measure on $\mathbb{T}$, $n$ times with itself. 
	Let us fix $N\geq 2$ and $m\in\mathbb{Z}_+$. We define
	$$L(z_1,...,z_N)=\frac{z_1+...+z_N}{\sqrt{N}},$$ and for $z^{(j)}=(z_{(j-1)N+1},...,z_{jN}),$  $$F(z_1,...,z_{mN})=\prod_{j=1}^m L(z^{(j)}).$$
	For $t>0$ and $0\leq\eta\leq 1/6$, we define $$F_t=(1+tF)^2-\eta t^2F^2, \quad \phi_t=1-t^2+tF.$$  
	 Since $$\int_{\mathbb{T}^{mN}}|F|^2dm_{mN}=\prod_{j=1}^m\int_{\mathbb{T}^{N}}|L|^2dm_{N}=1,$$ and $$\int_{\mathbb{T}^{mN}}F|F|^2dm_{mN}=0,$$ we get that
	$$\langle F_t,\phi_t\rangle =1+t^2.$$
	Now, let us define, for $S\subset \{1,...,m\}$ 
	$$e_S=\prod_{j\in S}L(z^{(j)}).$$
	Then the functions $e_S$ are orthonormal, and $\Ha_F e_S= e_{S^c}$. $\Ha_{\phi_t}$ acts on $\spans\{e_S,S\subset \{1,...,m\}\}$ since $H_L$ acts on $\spans\{1,L\}$, and its matrix decomposes into blocks corresponding to the pairs $\{e_S,e_{S^c}\}$. For $S=\emptyset$ the block is \[
	\left(
	\begin{array}{cc}
		1-t^2 & t \\
		t & 0
	\end{array}
	\right)
	\]
	and for every other pair ($2^{m-1}-1$ in total), the block is \[
	\left(
	\begin{array}{cc}
		0& t \\
		t & 0
	\end{array}
	\right).
	\]
	Therefore, we have that
	$$\|\Ha_{\phi_t}\|_{S^p}^p=(1+t^{2p})+(2^{m-1}-1)(2t^p)=1+t^{2p}+(2^{m}-2)t^p.$$
	Thus, equation \eqref{eq2} gives
	\begin{equation}\label{eq3} \|F_t\|_{L^1(\mathbb{T}^{mN})}\geq \dfrac{1+t^2}{(1+t^{2p}+(2^m-2)t^p)^{1/p}}.
	\end{equation}
	Let $$D_\eta(u)=1+u^2-\int_{\mathbb{T}}|(1+u\zeta)^2-\eta u^2\zeta^2|dm_1(\zeta).$$
	Now we can write
	$$\|F_t\|_{L^1}=1+t^2-\int_{\mathbb{T}^{mN}}D_\eta(t|F(z)|)dm_{mN}(z),$$
	and thus, substituting the norm into equation \eqref{eq3} and integrating the inequality against $t^{-\kappa-1}dt$, for $2< \kappa< 4$ (the choice of the range $(2,4)$ was made to make sure that both the left-hand side and the right-hand side of inequality \eqref{eq5} are finite), 
	we get that
	\begin{align}\label{eq4}
		\int_0^\infty \int_{\mathbb{T}^{mN}}\frac{D_\eta(t|F(z)|)}{t^{\kappa+1}}dm_{mN}(z) dt & \nonumber \\ \leq& \int_0^\infty \frac{(1+t^2)(1-(1+t^{2p}+(2^m-2)t^p)^{-1/p})}{t^{\kappa+1}}dt.
	\end{align}
	Substituting $t=2^{-m/p}x$, the right-hand side becomes
	$$2^{m\kappa/p}\int_0^\infty \frac{(1+2^{-2m/p}x^2)(1-(1+2^{-2m}x^{2p}+(2^m-2)2^{-m}x^p)^{-1/p})}{x^{\kappa+1}}dx.$$
	For $x<1$, since for $y>0$, $1-(1+y)^{-1/p}\leq y$, we have that
	$$1-(1+2^{-2m}x^{2p}+(2^m-2)2^{-m}x^p)^{-1/p}\leq (2^{-2m}x^{2p}+(2^m-2)2^{-m}x^p)\leq2x^p,$$
	thus for
	$$I_1:=\int_0^1 \frac{(1+2^{-2m/p}x^2)(1-(1+2^{-2m}x^{2p}+(2^m-2)2^{-m}x^p)^{-1/p})}{x^{\kappa+1}}dx,$$ we have $$I_1 \leq 4\int_0^1 x^{p-\kappa-1}dx=\frac{4}{p-\kappa}.$$
	For $x>1$ we use the simpler fact that
	$$1-(1+2^{-2m}x^{2p}+(2^m-2)2^{-m}x^p)^{-1/p}\leq 1,$$ thus getting for
	$$I_2:=\int_1^\infty \frac{(1+2^{-2m/p}x^2)(1-(1+2^{-2m}x^{2p}+(2^m-2)2^{-m}x^p)^{-1/p})}{x^{\kappa+1}}dx,$$ the upper bound
	$$I_2\leq \int_1^\infty \frac{1+x^2}{x^{\kappa+1}}dx=\frac{1}{\kappa}+\frac{1}{\kappa-2}.$$
	
    Thus equation \eqref{eq4} becomes
	\begin{equation}\label{eq5}	\int_0^\infty \int_{\mathbb{T}^{mN}}\frac{D_\eta(t|F(z)|)}{t^{\kappa+1}}dm_{mN}(z) dt \leq C_{p,\kappa}2^{m\kappa/p}.
	\end{equation}  
	We now distinguish between the cases $u\leq1/\sqrt{2}$ and $u\geq1/\sqrt{2}$. This distinction will ensure integrability of $D_\eta(u)/u^{\kappa+1}$ near $0$ and $\infty$, while the specific choice $1/\sqrt{2}$ will guarantee that the coefficient of $\eta$ in inequality \eqref{eq11} is positive.
	
	 To verify this, for $u>1/\sqrt 2$ we can see that
	\begin{eqnarray}  
		|D_\eta(u)|&\leq & \left|1+u^2-\int_{\mathbb{T}}|(1+u\zeta)^2-\eta u^2\zeta^2|dm_1(\zeta)\right| \nonumber \\
		&= & \left|\int_{\mathbb{T}}|1+u\zeta|^2dm_1(\zeta)-\int_{\mathbb{T}}|(1+u\zeta)^2-\eta u^2\zeta^2|dm_1(\zeta)\right| \nonumber \\
		&\leq & \int_{\mathbb{T}}|\eta u^2\zeta^2|dm_1(\zeta)=\eta u^2. \nonumber 
	\end{eqnarray} 
	
	For $u\leq1/\sqrt 2$ we observe that $(1+z)^2-\eta z^2$ is a nonvanishing analytic function on $D(0,\frac{1}{1+\sqrt{\eta}})$, and thus has an analytic square root $g$. Since $$g(z)-1-z=\frac{-\eta z^2}{g(z)+1+z}$$ and $$|g(z)+1+z|=|1+z|\left| 1+\sqrt{1-\frac{\eta z^2}{(1+z)^2}}\right|\geq |1+z|\geq 1-1/\sqrt 2,\; |z|<1/\sqrt 2,$$ we get that
	\begin{eqnarray}  \label{eq6}
		|D_\eta(u)|&=& \left|1+u^2-\int_{\mathbb{T}}|(1+u\zeta)^2-\eta u^2\zeta^2|dm_1(\zeta)\right| \nonumber \\
		&\leq& \int_{\mathbb{T}}|g(u\zeta)-1-u\zeta|^2dm_1(\zeta)\leq12 \eta^2 u^4.
	\end{eqnarray} 
	Also, by Fubini's theorem and the change of variables $u=t|F(z)|$, the left-hand side of inequality \eqref{eq5} becomes
	$$\int_0^\infty \frac{D_{\eta}(u)}{u^{\kappa+1}}du\left(\int_{\mathbb{T}^N}|L|^\kappa dm_N\right)^m,$$ and thus 
	\begin{equation}\label{eq7}
		\int_0^\infty \frac{D_{\eta}(u)}{u^{\kappa+1}}du\left(\int_{\mathbb{T}^N}|L|^\kappa dm_N\right)^m\leq C_{p,\kappa}2^{m\kappa/p}. 
	\end{equation}
	Suppose that for each $\kappa$ we can find some $\eta$ small enough such that $$\int_0^\infty \frac{D_{\eta}(u)}{u^{\kappa+1}}du>0.$$ Then, taking $m$-th roots and letting $m\to \infty$ in inequality \eqref{eq7} we get that
	\begin{equation}\label{eq8}
		\int_{\mathbb{T}^N}|L|^\kappa dm_N\leq 2^{\kappa/p}. 
	\end{equation}
	Taking $\kappa\to 4$, this gives
	$$2^{4/p}\geq \int_{\mathbb{T}^N}|L|^4 dm_N=2-\frac{1}{N}.$$ Finally, taking $N\to \infty$ we get that $2\leq 2^{4/p}$, which is our contradiction, since $p>4$. To finish the proof, it suffices to prove that for every $\kappa$ we can find an $\eta\in (0,1/6)$ with $$I_\kappa(\eta):=\int_0^\infty \frac{D_{\eta}(u)}{u^{\kappa+1}}du>0.$$ Cauchy--Schwarz and the elementary inequality $\sqrt{x^2+y}\leq x+y/(2x)$ give
	\begin{eqnarray} \label{eq9}D_{\eta}(u)&\geq& 1+u^2-\sqrt{(1+(1+\sqrt{\eta})^2u^2)(1+(1-\sqrt{\eta})^2u^2)} \nonumber \\
		&=& 1+u^2-\sqrt{(1+u^2)^2+2\eta u^2(1-u^2)+\eta^2 u^4} \nonumber  \\
		&\geq& -\frac{2\eta u^2(1-u^2)+\eta^2 u^4}{2(1+u^2)}  
	\end{eqnarray} 
	
	When $u\leq 1/\sqrt2$, inequality \eqref{eq6} gives
	\begin{eqnarray}\label{eq10} 
		\int_0^{1/\sqrt2} \frac{D_{\eta}(u)}{u^{\kappa+1}}du\geq -12\eta^2 \int_0^{1/\sqrt2} u^{3-\kappa}du=\frac{-12\eta^2}{2^{(4-\kappa)/2}(4-\kappa)}\geq \frac{-12\eta^2}{4-\kappa}.
	\end{eqnarray}
	For $u\geq 1/\sqrt 2$, inequality \eqref{eq9} and the observation that $$ \frac{u^2-1}{u^{\kappa-1}(1+u^2)}\geq  \frac{u^2-1}{u^3(1+u^2)},$$ give
	\begin{eqnarray}\label{eq11}
		\int_{1/\sqrt2}^\infty \frac{D_{\eta}(u)}{u^{\kappa+1}}du&\geq&\int_{1/\sqrt2}^\infty -\frac{2\eta u^{1-\kappa}(1-u^2)+\eta^2 u^{3-\kappa}}{2(1+u^2)}du \nonumber \\
		&=&\eta\int_{1/\sqrt2}^\infty \frac{ u^{1-\kappa}(u^2-1)}{1+u^2}du -C_\kappa \eta^2 \nonumber \\
		&\geq&\eta\int_{1/\sqrt2}^\infty \frac{u^2-1}{u^3(1+u^2)}du -C_\kappa \eta^2 \nonumber \\
		&=&\frac{\eta}{2}\int_{0}^2 \frac{1-x}{x+1}dx  -C_\kappa \eta^2 \nonumber \\
		&=&\eta(\log 3-1)  -C_\kappa \eta^2.
	\end{eqnarray}
	Here, for the second equality, we used the change of variables $x=u^{-2}$. Combining inequalities \eqref{eq10} and \eqref{eq11} we get that
	$$I_\kappa(\eta)\geq \eta(\log 3-1)-C_\kappa' \eta^2.$$ Thus we can choose $0<\eta<\min\left\{\frac16,\frac{\log3-1}{C_\kappa'}\right\}$.

	\bibliographystyle{plain}
	\bibliography{Torus}
	
\end{document}